\documentclass[12pt,a4paper]{article}%
\usepackage{amsmath}
\usepackage{epsfig}%
\usepackage{amsfonts}%
\usepackage{amssymb}
\usepackage{enumitem}
\usepackage{float}
\RequirePackage{tikz}
\RequirePackage{cleveref}
\RequirePackage{url}
\RequirePackage{pgfplots}
\pgfplotsset{compat=newest}
\usepackage{algorithmicx}
\usepackage{algorithm}
\usepackage[noend]{algpseudocode}
\def\BState{\State\hskip-\ALG@thistlm}

\usepackage{color}

\usepackage[numbers]{natbib}
\usepackage[%
    font={small,sf},
    labelfont=bf,
    format=hang,    
    format=plain,
    margin=0pt,
    width=0.8\textwidth,
]{caption}
\tikzset{
    state/.style={
           rectangle,
           rounded corners,
           draw=black, very thick,
           minimum height=2em,
           inner sep=2pt,
           text centered,
           },
}
\usepackage{subcaption}

\usepackage{graphicx}
\usepackage[normalem]{ulem}
 \usepackage{authblk}
 
\newtheorem{theorem}{Theorem}[section]

\newtheorem{definition}[theorem]{Definition}

\crefname{lemma}{Lemma}{Lemmas}

\newtheorem{remark}{Remark}
\crefname{remark}{Remark}{Remarks}

\makeatletter\@addtoreset{equation}{section}\makeatother
\makeatletter\@addtoreset{figure}{section}\makeatother
\makeatletter\@addtoreset{table}{section}\makeatother
\newenvironment{rcases}
  {\left.\begin{aligned}}
  {\end{aligned}\right\rbrace}

\begin{document}


\title{Solution decomposition for the nonlinear Poisson-Boltzmann equation using the range-separated 
tensor format}
 
\author[1,3]{Cleophas Kweyu \thanks{\tt kweyu@mpi-magdeburg.mpg.de}}
\author[1,2]{Venera Khoromskaia \thanks{\tt vekh@mis.mpg.de}}
\author[2]{Boris Khoromskij \thanks{\tt bokh@mis.mpg.de}}
\author[1]{Matthias Stein \thanks{\tt matthias.stein@mpi-magdeburg.mpg.de}}
\author[1]{Peter Benner \thanks{\tt benner@mpi-magdeburg.mpg.de}}
\affil[1]{Max Planck Institute for Dynamics of Complex Technical Systems, Sandtorstr.~1, D-39106 Magdeburg, 
Germany}
\affil[2]{Max Planck Institute for Mathematics in the Sciences, Inselstr.~22-26, D-04103 Leipzig, Germany} 
\affil[3]{Moi University, Department of Mathematics and Physics, P.O. Box 3900-30100, Eldoret, Kenya}


\maketitle

\begin{abstract}
The Poisson-Boltzmann equation (PBE) is an implicit solvent continuum model for 
calculating the electrostatic potential and energies of charged biomolecules in 
ionic solutions. However, its numerical solution poses a significant challenge due 
strong singularities and nonlinearity caused by the singular source terms and the 
exponential nonlinear terms, respectively. An efficient method for the treatment of 
singularities in the linear PBE which was introduced in \cite{BeKKKS:18}, that is based 
on the range-separated (RS) tensor decomposition \cite{BKK_RS:18} for both electrostatic 
potential and the discretized Dirac delta distribution \cite{khor-DiracRS:2018}. 
In this paper, we extend this regularization method to the  nonlinear PBE. 
Similar to \cite{BeKKKS:18} we apply the PBE only to the regular part of the solution corresponding to the 
modified right-hand side via extraction of the long-range part in the discretized Dirac delta distribution. 
The total electrostatic potential is obtained by adding the long-range solution to the directly precomputed 
short-range potential. The main computational benefit of the approach is the automatic preservation of the 
continuity in the Cauchy data on the solute-solvent interface. The boundary conditions are also obtained 
from the long-range component of the precomputed canonical tensor representation of the Newton kernel. In 
the numerical experiments, we illustrate the accuracy of the nonlinear regularized PBE (NRPBE) over the 
classical variant.
\end{abstract}



\section{Introduction}
\label{sec:Intro}

Biochemical processes are occurring between macromolecules such as proteins and nucleic acids in 
solution at a physiological salt concentration. The resultant electrostatic interactions are 
highly relevant for an understanding of biological functions and structures of biomolecules, enzyme catalysis, molecular recognition, and 
biomolecular encounter or association rates \cite{Fogolari2002,Neves-Petersen2003,Stein:2007,Stein:2010}. 
Efficient modeling of these interactions remains a great challenge in computational biology because 
of the complexity of biomolecular systems which are dominated by the effects of solvation on 
biomolecular processes and by the long-range intermolecular interactions 
\cite{DesHolm:98,LiStCaMaMe:13,ren2012}. 

There are two main types of models which can be used to model electrostatic 
interactions in ionic solutions. The explicit approaches which treat both the solute and solvent in 
atomic detail, are generally computationally demanding. This is because they require substantial 
sampling and equilibration in order to converge the properties of interest in an ensemble 
of solute and solvent \cite{ren2012,Jurrus2018}. On the other hand, continuum or 
implicit approaches treat the solvent molecules as a continuum, by integrating out 
non-relevant degrees of freedom in order to circumvent 
the need for sampling and equilibration \cite{Bashford:2000,ren2012,Jurrus2018}.

\begin{figure}[b!]
  \centering
    \includegraphics[width=12.0cm]{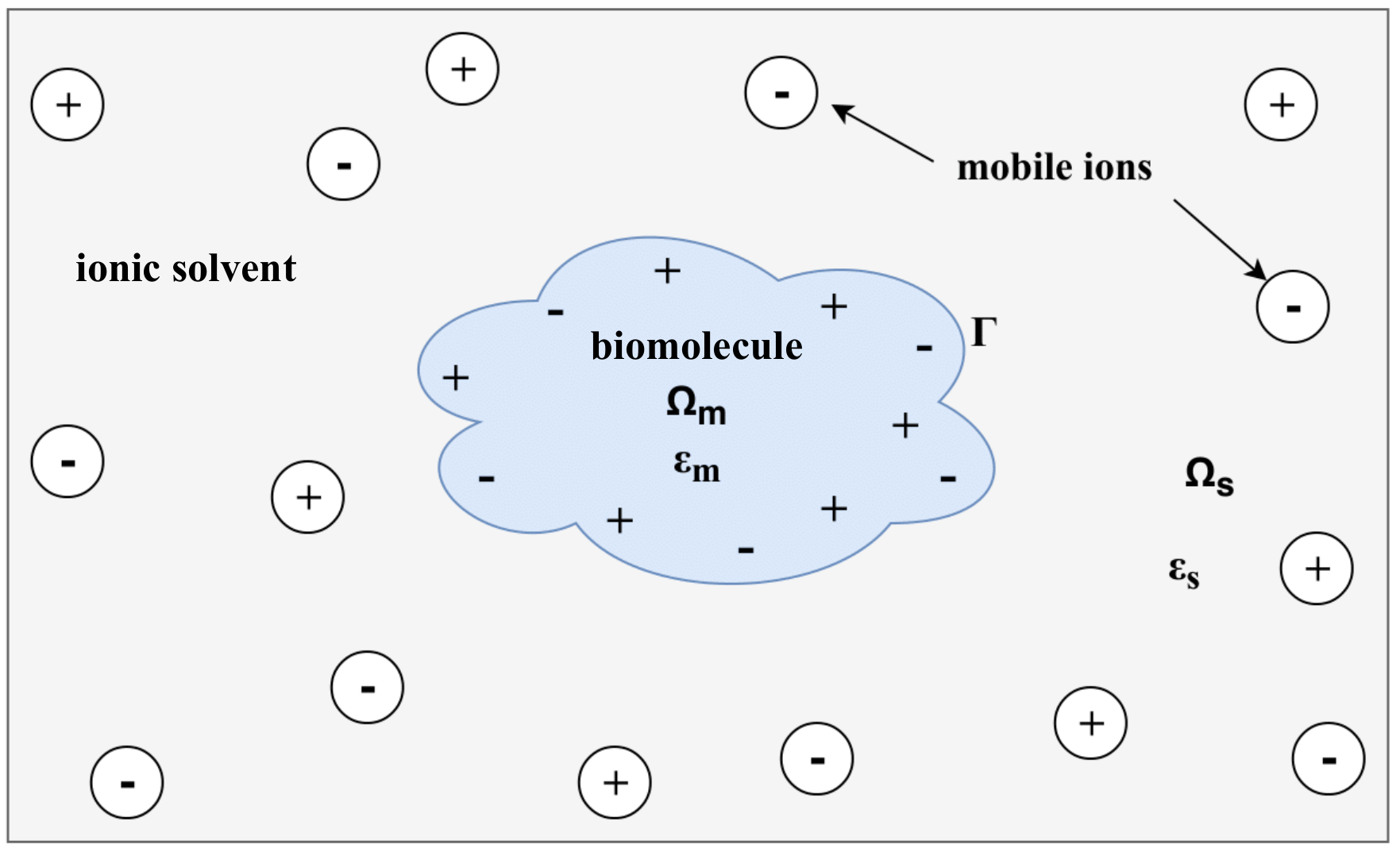}
    \caption{\label{fig:Biomolecular_system}
    Representation of a biomolecule with partial charges and an internal 
    dielectric constant ${\epsilon}_m$ in a solvent with dielectric constant 
    ${\epsilon}_s$ of mobile ions. Figure courtesy \cite{KwBFM2022}.}
\end{figure}

 There exists a number of implicit solvation approaches for biomolecules \cite{Bashford:2000,Barone:1997,
ren2012}, but the most popular is based on the Poisson-Boltzmann equation (PBE), which was extensively 
analyzed, for example, in \cite{Holst94}. The PBE is used for calculating the electrostatic potential and 
energies of charged or partially charged biomolecules in a physiological environment. 
We present the PBE model in \Cref{sec:PBE_theory}. 

It is impossible to obtain analytic solutions of the PBE for biomolecules with complex geometries and 
highly singular charge distributions \cite{Holst94, Dong2008}. The numerical solution of the PBE was 
pioneered by Warwicker and Watson in 1982 \cite{Warwicker1982}, where the electrostatic potential was 
computed at the active site of an enzyme using the finite difference method (FDM). Besides the FDM 
\cite{Baker2001, Wang2010}, other numerical techniques such as the finite element methods (FEM) 
\cite{Baker2001, Holst2000} and the boundary element methods (BEM) \cite{Boschitsch2004, Zhou1993} have 
hitherto successfully been used to solve the PBE, see \cite{Holst:2008} for a thorough review. However, the 
numerical solution of the PBE is faced with a number of challenges. The most significant are the strong 
charge singularities caused by the singular source terms (Dirac delta distribution), the nonlinearity 
caused by the exponential nonlinear terms, the unbounded domain due to slow polynomial decay of the 
potential with respect to distance, and by imposing the correct jump or interface conditions 
\cite{Xie:14, Mirzadeh:13}. 

The presence of a highly singular right-hand side of (\ref{eqn:PBE}) which is described by a sum of Dirac 
delta distributions introduces significant errors in the numerical solution of the PBE. To overcome this 
problem, the PBE theory has recently received a major boost by the introduction of solution decomposition 
(regularization) techniques which have been developed, for example, in 
\cite{Xie:14,Mirzadeh:13,Chen:07,Chern:2003}, see the discussion in \Cref{sec:soln_decomposition}. 
The idea behind these regularization techniques is the avoidance of building numerical 
approximations corresponding to the Dirac delta distributions by treating the biomolecular system (see 
\Cref{fig:Biomolecular_system}), as an interface problem. This is coupled with the advantage that analytical 
expansions in the molecular sub-region are possible, by the Newton kernel.

In this paper, for resolving the problem of strong singularities, we apply the method introduced recently 
in \cite{BeKKKS:18}, for the computation of the free-space electrostatic potential of a linear PBE and 
Poisson equation. For this purpose the range-separated canonical tensor format was applied, which was 
introduced and analyzed in \cite{BKK_RS:18, BKK_RS:16}. We extend the results of \cite{BeKKKS:18} to the 
case of   nonlinear PBE and compare the method numerically for a number of biomolecules. Similar to 
\cite{BeKKKS:18}, we apply the PBE only to the regular part of the solution corresponding to the 
modified right-hand side via extraction of the long-range part in the discretized Dirac delta 
distribution \cite{khor-DiracRS:2018}. Other numerical methods for the efficient treatment of the 
long-range part in the multi-particle electrostatic potential have been considered in \cite{Grigori:22}. 
 
The RS tensor formats can be gainfully applied to computational problems which include functions with 
multiple local singularities or cusps, Green kernels with intrinsic non-local behavior, and in various 
approximation problems which are generated by radial basis functions. The grid-based canonical tensor 
representation for the Newton kernel was developed in \cite{BeHaKh:08} and then applied in 
tensor-based electronic structure calculations \cite{khor-ml-2009,Khor_bookQC_2018}. Tensor numerical 
techniques for super-fast computation of the collective electrostatic potentials of large finite lattice 
clusters have been previously introduced in \cite{VeBoKh:Ewald:14}.

The splitting technique employed in this paper is based on the RS tensor decomposition of the discretized 
Dirac delta distribution \cite{khor-DiracRS:2018}, which allows avoiding the nontrivial matrix 
reconstruction as in (\ref{eqn:PBE_solndecomp}) and in \cite{Xie:14}. The only requirement in this approach 
is a simple modification of the singular charge density of the PBE in the molecular region $\Omega_m$, 
which does not change the FEM/FDM system matrix. The singular component in the total potential is recovered 
explicitly by the short-range component in the RS tensor splitting of the Newton potential. The main 
computational benefits of this approach are the localization of the modified singular charge density 
within the molecular region while automatically maintaining the continuity in the Cauchy data on the 
interface. Furthermore, this computational scheme only includes solving a single system of FEM/FDM 
equations for the regularized (or long-range) component of the decomposed potential.

The remainder of this paper is structured as follows. \Cref{sec:RS_survey} describes the 
basic rank-structured tensor formats and the short description of the range-separated tensor format 
\cite{BKK_RS:18, BKK_RS:16} for representation of the electrostatic potential of multiparticle systems. 
\Cref{sec:soln_decomposition} provides insights into the existing solution decomposition techniques for 
the PBE model. \Cref{sec:RS_2_PBE} explains how the application of the RS tensor format leads to the new regularization 
scheme for solving the PBE. \Cref{sec:RPBE_numerical_approach} presents the numerical
approach of solving the NRPBE. Finally, \Cref{sec:Numer_Tests} presents the numerical tests illustrating the 
benefits of the proposed method and comparisons with the solutions obtained by the 
standard FEM/FDM-based 
PBE solvers.

\section{The Poisson-Boltzmann equation theory}
\label{sec:PBE_theory}

The PBE is a nonlinear elliptic partial differential equation (PDE) which computes a global solution for 
the electrostatic potential within the biomolecule ($\Omega_m$) and in the surrounding ionic solution 
($\Omega_s$), see \Cref{fig:Biomolecular_system} for an illustration of the two regions. For a monovalent 
electrolyte (i.e., $1:1$ ion ratio), the dimensionless PBE is given by 
\begin{equation}\label{eqn:PBE}
    -\nabla\cdot(\epsilon(\bar{x})\nabla u(\bar{x})) + \bar{\kappa}^2(\bar{x})\sinh(u(\bar{x})) = 
      \sum_{i=1}^{N_m}q_i\delta(\bar{x}-\bar{x}_i),\quad \Omega \in 
    \mathbb{R}^3,
    \end{equation}
    subject to 
    \begin{equation}\label{eq:DH_solution}
    u(\bar{x}) = \sum_{i=1}^{N_m}\frac{q_ie^{-\kappa(d-a_i)}}{4\pi\epsilon_s (1+\kappa a_i)d} 
    \quad \mbox{on} \, \, \partial{\Omega}, \quad d = \lVert \bar{x}-\bar{x}_i \rVert, 
    \quad \bar{x} = (x,y,z),     
    \end{equation}
where $u(\bar{x}) = {e_c\psi(\bar{x})}/{\kappa_B T}$ represents the dimensionless potential, $\psi(\bar{x})$ 
is the original electrostatic potential in centimeter-gram-second (cgs) units scaled to the thermal voltage 
$(\kappa_B T)/e_c$, $q_i = \frac{4\pi e_c^2}{\kappa_B T}z_i$, $N_m$ is the total number of partial point charges in the 
biomolecule, $\epsilon_s$ is the bulk solvent dielectric coefficient, and $a_i$ is the atomic radius of the mobile ions. 
Here, $\kappa_B T$, $\kappa_B$, $T$, $e_c$, and $z_i$ are the thermal energy, the Boltzmann constant, the absolute 
temperature, the electron charge, and the non-dimensional partial charge of each atom, respectively. The 
Debye-H\"uckel screening parameter, $\kappa^2 = {8\pi e_c^2 I}/{1000\epsilon_s \kappa_BT}$, describes the ion 
concentration and accessibility, and is a function of the ionic strength 
$I = 1/2\sum_{j=1}^{N_{ions}}c_jz_j^2$, where $z_j$ and $c_j$ are charge and concentration of each ion. 
The sum of Dirac delta distributions, located at atomic centers $\bar{x}_i$, represents the molecular charge 
density. See \cite{Holst94, KwBFM2017} for more details concerning the PBE theory.

The dielectric coefficient $\epsilon(\bar{x})$ and kappa function $\bar{\kappa}^2(\bar{x})$ are piecewise 
constant functions given by
\begin{eqnarray}\label{eqn:diel_kappa_def}
 \epsilon(\bar{x}) =
  \begin{cases}
   \epsilon_m = 2 & \text{if } \bar{x} \in \Omega_m\\
   \epsilon_s \,\,= 78.54 & \text{if } \bar{x} \in \Omega_s
  \end{cases}, \quad \quad
 \bar{\kappa}(\bar{x}) =
  \begin{cases}
   0 & \text{if } \bar{x} \in \Omega_m \\
    \sqrt{\epsilon_s}\bar{\kappa} & \text{if } \bar{x} \in \Omega_s 
  \end{cases},
\end{eqnarray}
where $\Omega_m$ and $\Omega_s$ are the molecular and solvent regions, respectively, as shown in 
\Cref{fig:Biomolecular_system}. Details of regarding the PBE theory and the significance of (\ref{eqn:PBE}) 
in biomolecular modeling can be found in \cite{SharpHonig90,Holst94,KwBFM2017}.

The PBE in (\ref{eqn:PBE}) can be linearized for small electrostatic potentials relative to the thermal 
energy (i.e., $\psi(\bar{x}) \ll \kappa_BT$). Nevertheless, even when the linearization condition does not hold, the 
solution obtained from the linearized PBE (LPBE) is close to that of the nonlinear PBE \cite{Fogolari99}. 
The onset of substantial differences between the two models is attributed to the magnitude of the electric 
field, hence, of the charge density at the interface between the solute and the solvent \cite{Fogolari99}. 
The LPBE is given by
\begin{equation}\label{eqn:LPBE}
   -\nabla\cdot(\epsilon(\bar{x})\nabla u(\bar{x})) + \bar{k}^2(\bar{x})u(\bar{x}) 
				      = \sum_{i=1}^{N_m}q_i\delta(\bar{x}-\bar{x}_i).
\end{equation}

The electrostatic potential can be used in a variety of applications, a few of which we highlight here. 
First, the surface potential, (i.e., the electrostatic potential on the biomolecular surface), can be 
used to obtain insights into possible binding sites for other molecules. Secondly, it can be used to compare 
the interaction properties of related proteins by calculating similarity indices \cite{Wade_etal:2001}. 
Finally, the electric field, which is the derivative of the potential around the solute, may be essential for 
obtaining the rates of molecular recognition and encounter \cite{Fogolari2002,Jurrus2018}.

\section{Rank-structured tensor representation of electrostatic potentials}
\label{sec:RS_survey}

\subsection{Sketch of basic tensor formats}
 \label{ssec:Tensor_formats}

Here, we recall the rank-structured tensor formats and briefly describe the range-separated tensor format 
introduced in \cite{BKK_RS:18,BKK_RS:16} for tensor-based representation of multiparticle long-range 
potentials. Rank-structured tensor techniques have recently gained popularity in scientific computing due 
to their inherent property of reducing the grid-based solution of multidimensional problems arising in 
large-scale electronic and molecular structure calculations to essentially 1D computations 
\cite{khor-ml-2009,KhKhFl_Hart:09}. In this concern, the so-called reduced higher order singular value 
decomposition (RHOSVD) introduced in \cite{khor-ml-2009} is one of the salient ingredients in the 
development of tensor methods in quantum chemistry, see details in \cite{Khor_bookQC_2018} and references 
therein.
 
A tensor of order $d$ is defined as a real multidimensional array over a $d$-tuple index set
 \begin{equation}
     {\bf A} = [a_{i_1, \ldots, i_d}] \equiv [a(i_1, \ldots, i_d)] \in \mathbb{R}^{n_1 
     \times \cdots \times n_d},
    \end{equation}
  with multi-index notation $i= (i_1, \ldots, i_d)$, $i_{\ell} \in I_{\ell} := \{1, \dots, 
  n_{\ell}\}$. It is considered as an element of a linear vector space $\mathbb{R}^{n_1\times \cdots 
  \times n_d}$ equipped with the Euclidean scalar product $\langle \cdot,\cdot \rangle : \mathbb{V}_n 
\times \mathbb{V}_n \rightarrow \mathbb{R}$, defined as
\begin{equation}
 \langle {\bf A},{\bf B} \rangle := \sum\limits_{(i_1, \ldots, i_d) \in I} a_{i_1, \ldots, i_d} 
 b_{i_1, \ldots, i_d} \quad \mbox{for}\quad {\bf A}, \,{\bf B} \, \in \mathbb{V}_n.
\end{equation}
The storage size scales exponentially in the dimension $d$, i.e., $n^d$, resulting in the so-called ``curse of 
dimensionality''. To get rid of the exponential scaling in storage and the consequent drawbacks,
one can apply the rank-structured separable approximations of multidimensional tensors. 
The simplest separable tensor is given by a rank-1 canonical tensor (i.e., tensor/outer product of 
 vectors in $d$ dimensions)
 \begin{equation}\label{eqn:rank1_canon}
  {\bf U} = {\bf u}^{(1)} \otimes \cdots \otimes {\bf u}^{(d)}  
  \in \mathbb{R}^{n_1 \times \cdots \times n_d},
\end{equation}
with entries computed as $u_{i_1, \ldots, i_d} = u_{i_1}^{(1)}  \cdots u_{i_1}^{(d)}$, which requires only 
$(n_1 + \ldots + n_d)\ll n^d$ numbers to store it. If $n_{\ell} = n$, then the storage cost is $dn \ll n^d$. 
 
\begin{definition}
 The \textbf{$R$-term canonical tensor} format is defined by a finite sum of rank-1 tensors
 \begin{equation}\label{eqn:canonical}
     {\bf U}_R =  \sum_{k =1}^{R} \xi_k {\bf u}_k^{(1)}  \otimes \cdots \otimes {\bf u}_k^{(d)}, 
     \quad \xi_k \in \mathbb{R},
    \end{equation}
 where ${\bf u}_k^{(\ell)} \in \mathbb{R}^{n_{\ell}}$ are normalized vectors, and $R \in \mathbb{R}_+$ is 
 the canonical rank. 
 \end{definition}
 The storage cost for this tensor format is bounded by $dRn$. For $k=3$, for example, the entries of the 
 canonical tensor (\ref{eqn:canonical})  are computed as the sums of elementwise products,  
 \begin{equation}
  u_{i_1,i_2,i_3} = \sum_{k=1}^R \xi_k u_{i_1,k}^{(1)} \cdot u_{i_2,k}^{(2)} \cdot u_{i_3,k}^{(3)}.
 \end{equation}
  
\begin{definition}
The \textbf{rank-${\bf r}$ orthogonal Tucker format} for a tensor ${\bf V}$ is
    \begin{equation}\label{eqn:Tucker}
 {\bf V} =  \sum_{\nu_1 =1}^{r_1} \cdots \sum_{\nu_d =1}^{r_d} \beta_{\nu_1, \ldots, \nu_d} 
 {\bf v}_{\nu_1}^{(1)}  \otimes \cdots \otimes {\bf v}_{\nu_d}^{(d)}  \equiv  
 \boldsymbol{\beta} \times_1 V^{(1)} \times_2 V^{(2)}\ldots \times_d V^{(d)},
 \end{equation}
 where $\{ {\bf v}_{\nu_{\ell}}^{(\ell)}\}_{\nu_{\ell} = 1}^{r_{\ell}} \in \mathbb{R}^{n_{\ell}}$ 
 is the set of orthonormal vectors for $\ell = 1, \ldots, d$.
$\times_{\ell}$ denotes  the contraction along the mode $\ell$ with the orthogonal matrices 
 $V^{(\ell)} = [{\bf v}_1^{(\ell)}, \ldots, {\bf v}_{r_{\ell}}^{(\ell)}] \in \mathbb{R}^{n_{\ell} 
 \times r_{\ell}}$.
$\boldsymbol{\beta} = \beta_{\nu_1, \ldots, \nu_d} \in \mathbb{R}^{r_1 \times \cdots r_d}$ is 
 the Tucker core tensor. 
\end{definition}
The storage cost is bounded by $drn + r^d$ with $r = |r| := \mbox{max}_{\ell} 
 r_{\ell}$.

 Rank-structured tensor approximations provide fast multilinear algebra with linear complexity scaling in 
 the dimension $d$ \cite{BKK_RS:18}. For instance, for the given canonical tensor representation 
 (\ref{eqn:canonical}), Hadamard products, the Euclidean scalar product, and $d$-dimensional 
 convolution can be computed by univariate tensor operations in 1D complexity \cite{KhKh:06}.

\subsection{Outline on the RS tensor format for numerical modeling of multiparticle systems}
 \label{ssec:Coulomb}

In what follows, first recall the canonical tensor representation of the non-local Newton kernel 
$1/\|\bar{x}\|$, $\bar{x} \in \mathbb{R}^3$, by using sinc-quadratures and Laplace transform introduced 
in \cite{BeHaKh:08}. The corresponding theoretical basis was developed in seminal papers 
\cite{HaKhtens:04I,khor-rstruct-2006} on low-rank tensor product approximation of multidimensional 
functions and operators. According to above papers, the Newton kernel is approximated in a computational 
domain $\Omega=[-b,b]^3$, using the uniform $n\times n\times n$ 3D Cartesian grid. Then, using the Laplace 
transform and sinc-quadrature approximation, this discretized potential is approximated by a canonical 
rank $R$ tensor,
\begin{equation} \label{eqn:canon_repr}
\mathbf{P} \approx  
 \sum\limits_{k=1}^{R} {\bf p}^{(1)}_k \otimes {\bf p}^{(2)}_k \otimes {\bf p}^{(3)}_k
\in \mathbb{R}^{n^{\otimes 3}},
\end{equation}
with vectors ${\bf p}^{(\ell)}_k \in \mathbb{R}^n$, and the accuracy of this approximation decays 
exponentially fast in the rank parameter $R$.
 
The canonical tensor representation of the Newton kernel was first applied in rank-structured grid-based 
calculations of the multidimensional operators in electronic structure calculations, 
\cite{khor-ml-2009,KhKhAn:12}, where it manifested its high accuracy compared with analytical based 
computational methods. 

In \cite{VeBoKh:Ewald:14}, the canonical tensor representation was applied in modeling 
of the electrostatic potentials in finite rectangular three-dimensional lattices, where 
it was proven that the rank of the collective long-range electrostatic potentials 
of large 3D lattices remains as small as that of a canonical tensor for a single Newton kernel. 
For lattices with defects and impurities it is higher by a small constant \cite{Khor_bookQC_2018}.
    
For modeling the electrostatic interaction potential in large molecular systems of general type, the 
range-separated tensor format \cite{BKK_RS:18} is based on additive decomposition of the reference canonical 
tensor $\textbf{P}_R$ 
 \[
  \mathbf{P}_R = \mathbf{P}_{R_s} + \mathbf{P}_{R_l},
\]
with
\begin{equation} \label{eqn:Split_Tens}
    \mathbf{P}_{R_s} =
\sum\limits_{k\in {\cal K}_s} {\bf p}^{(1)}_k \otimes {\bf p}^{(2)}_k \otimes {\bf p}^{(3)}_k,
\quad \mathbf{P}_{R_l} =
\sum\limits_{k\in {\cal K}_l} {\bf p}^{(1)}_k \otimes {\bf p}^{(2)}_k \otimes {\bf p}^{(3)}_k.
\end{equation}

Here, ${\cal K}_l := \{k|k = 0,1, \ldots, R_l\}$ and ${\cal K}_s := \{k|k = R_l+1, \ldots, M\}$ are the 
sets of indices for the long- and short-range canonical vectors determined depending on the claimed size 
of effective support of the short-range part $\mathbf{P}_{R_s}$. 
  
The total electrostatic potential   is represented by a projected tensor 
${\bf P}_0\in \mathbb{R}^{n \times n \times n}$ that can be constructed by a direct sum of 
shift-and-windowing transforms of the reference tensor $\widetilde{\bf P}_R$, defined in the twice larger 
domain $\widetilde{\Omega}_n$ (see \cite{VeBoKh:Ewald:14} for more details),
\begin{equation}\label{eqn:Total_Sum}
 {\bf P}_0 = \sum_{\nu=1}^{N} {z_\nu}\, {\cal W}_\nu (\widetilde{\bf P}_R)=
 \sum_{\nu=1}^{N} {z_\nu} \, {\cal W}_\nu (\widetilde{\mathbf{P}}_{R_s} + \widetilde{\mathbf{P}}_{R_l})
 =: {\bf P}_s + {\bf P}_l.
\end{equation}
The shift-and-windowing transform ${\cal W}_\nu$ maps a reference tensor 
$\widetilde{\bf P}_R\in \mathbb{R}^{2n \times 2n \times 2n}$ onto its sub-tensor 
of smaller size $n \times n \times n$, obtained by first shifting the center of
the reference tensor $\widetilde{\bf P}_R$ to the grid-point $x_\nu$ and then restricting 
(windowing) the result onto the computational grid $\Omega_n$.

It was proven in \cite{BKK_RS:18} that the Tucker and canonical rank parameters of the "long-range part" 
in the tensor ${\bf P}_0$, defined by
\begin{equation}\label{eqn:Long-Range_Sum} 
 {\bf P}_l = \sum_{\nu=1}^{N} {z_\nu} \, {\cal W}_\nu (\widetilde{\mathbf{P}}_{R_l})=
 \sum_{\nu=1}^{N} {z_\nu} \, {\cal W}_\nu 
 (\sum\limits_{k\in {\cal K}_l} \widetilde{\bf p}^{(1)}_k \otimes \widetilde{\bf p}^{(2)}_k 
 \otimes \widetilde{\bf p}^{(3)}_k)
 \end{equation}
 remain almost uniformly bounded in the number of particles,
 \[
  \mbox{rank}({\bf P}_l)\leq C \log^{3/2} N.
 \]
The rank reduction algorithm is accomplished by the canonical-to-Tucker (C2T) transform through 
the reduced higher order singular value decomposition (RHOSVD) \cite{khor-ml-2009} with a 
subsequent Tucker-to-canonical (T2C) decomposition (see \cite{Khor_bookQC_2018} and references therein).

In turn, the tensor representation of the sum of short-range parts is considered as a sum of cumulative 
tensors of small support characterized by the list of the 3D potentials coordinates and weights. The 
total tensor is then represented in the range-separated tensor format \cite{BKK_RS:18}. 
Here, we recall a slightly simplified definition of the RS tensor format.
\begin{definition}\label{Def:RS-Can_format} (RS-canonical tensors \cite{BKK_RS:18}). 
Given a reference tensor ${\bf A}_0$ such that $\mbox{rank}({\bf A}_0)\leq R_0$,
the separation parameter $\gamma \in \mathbb{N}$ and a set of points $x_\nu \in \mathbb{R}^{d}$,
$\nu=1,\ldots,N$,
the \textbf{RS-canonical tensor format} specifies the class of $d$-tensors 
${\bf A}  \in \mathbb{R}^{n_1\times \cdots \times n_d}$
which can be represented as a sum of a rank-${R}_L$ canonical tensor  
\begin{equation}\label{eq:LR_tensor_sum}
{\bf A}_{R_L} = {\sum}_{k =1}^{R_L} \xi_k {\bf a}_k^{(1)} \otimes \cdots \otimes {\bf a}_k^{(d)}
\in \mathbb{R}^{n_1\times ... \times n_d}
 \end{equation}
and a cumulated canonical tensor 
 \begin{equation}\label{eq:CCT_tensor}
  \widehat{\bf A}_S={\sum}_{\nu =1}^{N} c_\nu {\bf A}_\nu ,  
 \end{equation}
generated by replication of the reference tensor ${\bf A}_0$ to the points $x_\nu$.
Then the RS canonical tensor is represented in the form
\begin{equation}\label{eqn:RS_Can}
 {\bf A} =  {\bf A}_{R_L} + \widehat{\bf A}_S=
 {\sum}_{k =1}^{R_L} \xi_k {\bf a}_k^{(1)}  \otimes \cdots \otimes {\bf a}_k^{(d)} +
 {\sum}_{\nu =1}^{N} c_\nu {\bf A}_\nu, 
\end{equation}
where $\mbox{diam}(\mbox{supp}\,{\bf A}_0)\leq 2 \gamma$ in the index size.
\end{definition}
The storage size for  the  RS-canonical tensor ${\bf A}$ in (\ref{eqn:RS_Can}) 
is estimated by (\cite{BKK_RS:18}, Lemma 3.9),
$$
\mbox{stor}({\bf A})\leq d R n + (d+1)N + d R_0 \gamma.
$$

Notice that the RS tensor decomposition of the collective electrostatic potential $\mathbf{P}_{0}$
can be obtained by setting  ${\bf A}_0=\mathbf{P}_s$ and 
${\bf A}_{R_L} = {\bf P}_l$.

\section{Solution decomposition techniques for the PBE}
\label{sec:soln_decomposition} 

The presence of the highly singular right-hand side of (\ref{eqn:PBE}) implies that every singular charge 
$z_i$ in (\ref{eqn:PBE}), the electrostatic potential $u(\bar{x})$ exhibits degenerate behavior at each atomic 
position $\bar{x}_i$ in the molecular region $\Omega_m$. To overcome this difficulty, the PBE theory has 
recently received a major boost by the introduction of solution decomposition techniques 
which entail a coupling of two equations for the 
electrostatic potential in the molecular ($\Omega_m$) and solvent ($\Omega_s$) regions, through the 
boundary interface \cite{Chen:07, Chern:2003}. The equation inside $\Omega_m$ is simply the Poisson 
equation, due to the absence of ions, i.e.,  
\begin{equation}\label{eqn:PE}
 -\nabla \cdot(\epsilon_m\nabla u) = \sum_{i=1}^{N_m}q_i \delta(\bar{x}-\bar{x}_i)  
	\quad \mbox{in} \,\, \Omega_m,
\end{equation}
On the other hand, there is absence of atoms in $\Omega_s$. Therefore, the density is purely given by the 
Boltzmann distribution
\begin{equation}\label{eqn:Boltzmann_distn}
 -\nabla \cdot(\epsilon_s\nabla u) + \bar{\kappa}^2\sinh(u) = 0  \quad \mbox{in} \,\, \Omega_s.
\end{equation}
The two equations (\ref{eqn:PE}) and (\ref{eqn:Boltzmann_distn}) are coupled together through the 
interface boundary conditions
\begin{equation}\label{eqn:interface_condn}
 \left[ u\right]_{\Gamma} = 0, \quad \mbox{and} \quad \left[ \epsilon \frac{\partial u}{\partial 
 n_{\Gamma}}\right]_{\Gamma} = 0,
\end{equation}
where $\Gamma := \partial\Omega_m = \partial\Omega_s \cap \Omega_m$ and $\left[ f\right]_{\Gamma} = 
\lim\limits_{t \longrightarrow 0} f(\bar{x}+tn_{\Gamma}) - f(\bar{x}-tn_{\Gamma})$. Here, $n_{\Gamma}$ 
denotes the unit outward normal direction of the interface $\Gamma$. 

Next, we highlight one of the solution decomposition techniques for the PBE in \cite{Chen:07} which 
provides the motivation for the RS tensor format demonstrated in this paper. It is also implemented as an 
option for the PBE solution in the well-known adaptive Poisson-Boltzmann software (APBS) package using the 
FEM \cite{Bakersept2001}. To deal with the singular source term represented by the sum of Dirac delta 
distributions in the PBE, the unknown solution $u(\bar{x})$ is decomposed as an unknown smooth function 
$u^r(\bar{x})$ and a known singular function $G(\bar{x})$, i.e.,
\begin{equation}
 u(\bar{x}) = G(\bar{x}) + u^r(\bar{x}),
\end{equation} 
where 
\begin{equation}\label{eqn:analytic_PE}
G(\bar{x})= \sum_{i=1}^{N_m}\frac{q_i}{\epsilon_m}\frac{1}{\|\bar{x} -\bar{x}_i\|},
\end{equation} 
is a sum of the Newton kernels ($1/\|\bar{x}\|$), which solves the Poisson equation (\ref{eqn:PE}) in 
$\mathbb{R}^3$. Substitute the decomposition into (\ref{eqn:PBE}), to obtain
\begin{equation}\label{eqn:PBE_solndecomp}
\begin{rcases}
\begin{aligned}
-\nabla \cdot(\epsilon\nabla u^r) + \bar{\kappa}^2(\bar{x})\sinh(u^r +G) 
   &= \nabla \cdot((\epsilon-\epsilon_m)\nabla G),  & \, \mbox{in} \,\, \Omega \\
u^r &= g - G & \, \mbox{on} \,\, \partial{\Omega},
\end{aligned}
\end{rcases}
\end{equation}
where $g(\bar{x})$ is the boundary condition obtained from (\ref{eq:DH_solution}). The PBE in 
(\ref{eqn:PBE_solndecomp}) is referred to as the regularized PBE (RPBE) in \cite{Chen:07}. Notice that 
the singularities of the Dirac delta distribution are transferred to $G$, which is known analytically, 
therefore, building the numerical approximation to $G$ is circumvented. Consequently, the cutoff 
coefficients $\bar{\kappa}$ and $\epsilon-\epsilon_m$ are zero in $\Omega_m$, where the degenerate 
behaviour is exhibited at each $\bar{x}_i$. This allows the RPBE to be a mathematically well-defined 
equation for the regularized solution $u^r$. It is important to note that away from the ${\bar{x}_i}$, the 
function $G$ is smooth \cite{Chen:07}. 

The RPBE in (\ref{eqn:PBE_solndecomp}) can further be decomposed into the linear and nonlinear components, 
$u^r(\bar{x}) = u^l(\bar{x}) + u^n(\bar{x})$, where $u^l(\bar{x})$ satisfies,
\begin{equation}\label{eqn:PBE_solndecompb}
\begin{rcases}
\begin{aligned}
-\nabla \cdot(\epsilon\nabla u^l) &= \nabla \cdot((\epsilon-\epsilon_m)\nabla G),  
& \, \mbox{in} \,\, \Omega \\
u^l &= 0 & \, \mbox{on} \,\, \partial{\Omega},
\end{aligned}
\end{rcases}
\end{equation}
and $u^n(\bar{x})$ satisfies 
\begin{equation}\label{eqn:PBE_solndecompc}
\begin{rcases}
\begin{aligned}
-\nabla \cdot(\epsilon\nabla u^n) + \bar{\kappa}^2(x)\sinh(u^n + u^l + G) &= 0,  
& \, \mbox{in} \,\, \Omega \\
u^n &= g-G & \, \mbox{on} \,\, \partial{\Omega}.
\end{aligned}
\end{rcases}
\end{equation}

However, the following computational challenges are inherent in the aforementioned techniques. First, due 
to regularization splitting of the solution by using the kappa and dielectric coefficients as cutoff 
functions, discontinuities at the interface arise. Therefore, interface or jump conditions need to be 
incorporated to eliminate the solution discontinuity (e.g., Cauchy data) at the interface of complicated 
sub-domain shapes. Consequently, the long-range components of the free space potential are not completely 
decoupled from the short-range parts at each atomic radius, in the ``so-called'' singular function $G$, 
in the molecular domain $\Omega_m$. Secondly, the Dirichlet boundary conditions, for example, in 
(\ref{eq:DH_solution}) have to be specified using some analytical solution of the LPBE. Thirdly, in 
solution decomposition techniques, see, for instance, \cite{Xie:14}, multiple algebraic systems for the 
linear and nonlinear boundary value problems have to be solved, thereby increasing the computational 
costs. Thirdly, the system matrix is modified because of incorporating the interface conditions and also, 
for instance, the smooth function ($G$), in the Boltzmann distribution term in (\ref{eqn:PBE_solndecomp}).

In this paper, we present a new approach for the regularization of the PBE by using the RS canonical 
tensor format.  

\section{The regularization scheme for the PBE via RS tensor format} \label{sec:RS_2_PBE}
 
In this section, we extend the approach introduced in \cite{BeKKKS:18} for linear PBE to the nonlinear case.
We present a new regularization scheme for the nonlinear PBE which is based on the range-separated 
representation of the highly singular charge density, described by the Dirac delta distribution 
in the target PBE (\ref{eqn:PBE}) \cite{khor-DiracRS:2018}.  
Similar to \cite{BeKKKS:18} we modify the right-hand 
side of the nonlinear PBE (\ref{eqn:PBE}) in such a way that 
the short-range part in the solution $u$ can be pre-computed 
independently by the direct tensor decomposition of the free space potential, and the initial elliptic 
equation (or the nonlinear RPBE) applies only to the long-range component of the total potential. The latter is a 
smooth function, hence the FDM/FEM approximation error can be reduced dramatically even on relatively 
coarse grids in 3D.

\subsection{Regularization scheme for the nonlinear PBE (NPBE)} \label{ssec:nonlin_PBE}

To fix the idea, we first consider the weighted sum of interaction potentials in a large $N$-particle 
system, generated by the Newton kernel, $1/{\|\bar{x}\|}$, at each charge location $\bar{x}_i$, 
$\bar{x}\in \mathbb{R}^3$, i.e.,
\begin{equation}\label{eqn:analytic_PE2}
G(\bar{x})= \sum_{i=1}^{N_m}\frac{q_i}{\epsilon_m}\frac{1}{\|\bar{x} -\bar{x}_i\|}, 
\end{equation}  
We recall that the sum of Newton kernels for a multiparticle system discretized by the $R$-term 
sum of Gaussian type functions living on the $n^{\otimes 3}$ tensor grid $\Omega_n$ is represented 
by a sum of long-range tensors in (\ref{eqn:Long-Range_Sum}) and a cumulated canonical tensor in (\ref{eq:CCT_tensor}), 
respectively.

Since it is well known that (\ref{eqn:analytic_PE2}) solves the Poisson equation analytically, i.e.,
\begin{equation}\label{eqn:PE_substit}
 -\nabla \cdot(\epsilon_m\nabla G(\bar{x})) = \sum_{i=1}^{N_m}q_i \delta(\bar{x}-\bar{x}_i)  
 \quad \mbox{in} \,\, \mathbb{R}^3,
\end{equation}
we can leverage this property in order to derive a smooth (regularized) representation, $f_r$, of the 
Dirac delta distributions in the right-hand side of (\ref{eqn:PE_substit}).  Consider the RS tensor 
splitting of the multiparticle Newton potential into a sum of long-range tensors ${\bf P}_l$ in 
(\ref{eqn:Long-Range_Sum}) and a cumulated canonical tensor ${\bf P}_s$ in (\ref{eq:CCT_tensor}), i.e., 
\begin{equation}\label{eqn:Newt_splitting}
 G(\bar{x}) = {\bf P}_s(\bar{x}) + {\bf P}_l(\bar{x}).
\end{equation}

Substituting each of the components of (\ref{eqn:Newt_splitting}) into the discretized Poisson equation, 
we derive the respective components of the molecular charge density 
(or the collective Dirac delta distributions) as follows
\begin{equation}\label{eqn:Dirac_splitting}
 f^s:= -A_{\Delta} {\bf P}_s, \quad 
 \mbox{and} \quad f^l:= -A_{\Delta} {\bf P}_l,
\end{equation}
where $A_{\Delta}$ is the 3D finite difference Laplacian matrix defined on the uniform rectangular grid as
\begin{equation}\label{eqn:Lapl_Kron3}
A_{\Delta} = \Delta_{1} \otimes I_2\otimes I_3 + I_1 \otimes \Delta_{2} \otimes I_3 + 
I_1 \otimes I_2\otimes \Delta_{3},
\end{equation}
where $-\Delta_\ell = h_\ell^{-2} \mathrm{tridiag} \{ 1,-2,1 \} \in \mathbb{R}^{n_\ell \times n_\ell}$, 
$\ell=1,2,3$, denotes the discrete univariate Laplacian and $I_\ell$, $\ell=1,2,3$, is the identity matrix 
in each dimension. 
See \cite{BeKKKS:18,khor-DiracRS:2018} for more details.

\begin{figure}[t]
\captionsetup{width=\linewidth}
    \centering
    \begin{subfigure}[b]{0.45\textwidth}
        \centering
        \includegraphics[width=\textwidth]{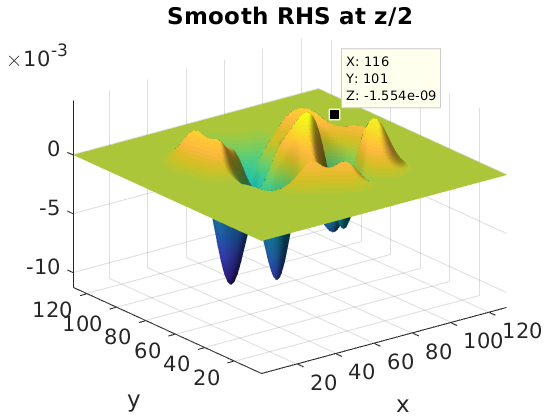}
        \caption{Long-range part of charge density.}
        \label{fig:Long-range_RHS}
    \end{subfigure}
    \hfill
    \begin{subfigure}[b]{0.45\textwidth}
        \centering
        \includegraphics[width=\textwidth]{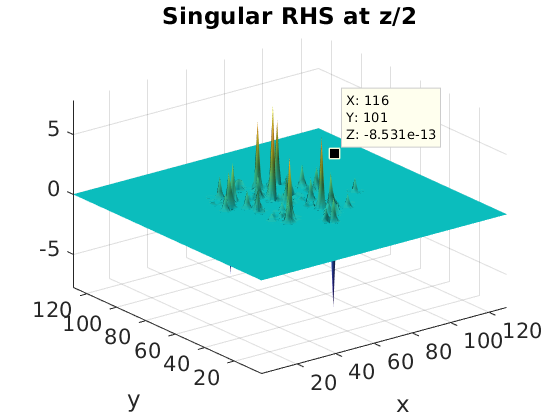}
        \caption{Short-range part of charge density.}
        \label{fig:Short-range_RHS}
    \end{subfigure}
    \caption{\label{fig:Long_short_RHS}
    The long- and short-range parts of the charge density for protein Fasciculin 1 on $129^{\otimes 3}$ 
    grid.}   
\end{figure}

\Cref{fig:Long_short_RHS} depicts the behaviour of the modified representations of both the smooth and 
singular components of the Dirac delta distributions using the formula in (\ref{eqn:Dirac_splitting}). The 
charge density data is obtained from protein Fasciculin 1, an anti-acetylcholinesterase toxin from green 
mamba snake venom \cite{DuMaBoFo:92}. Notice from the highlighted data cursors, that 
the effective supports of both functions 
are localized within the molecular region, with values dropping to zero outside this region. Furthermore, 
\Cref{fig:Long-range_RHS} represents the function $f^l$, which we utilize as the modified right-hand side
to derive a regularized PBE model (RPBE) in the next step.

The nonlinear regularized PBE (NRPBE) can now be derived as follows. First, the unknown solution 
(or target electrostatic potential) $u$ to the PBE (\ref{eqn:PBE}) can be decomposed as 
\[
u = u^s + u^r,
\]
where $u^s$ is the known singular function (or short-range component) and $u^l$ is the unknown long-range 
component to be determined. Therefore, the PBE (\ref{eqn:PBE}) can be rewritten as
\begin{equation}\label{eqn:PBE_splitting}
  \begin{rcases}
  \begin{aligned}
 -\nabla\cdot(\epsilon\nabla (u^s + u^r)) + \bar{\kappa}^2\sinh(u^s + u^r) &= f^s + f^l  
 \quad \mbox{in} \,\, \mathbb{R}^3,\\
 u^r &= g, \quad \mbox{on} \quad \partial\Omega,
  \end{aligned}
  \end{rcases}
\end{equation}
where the right-hand side of (\ref{eqn:PBE}) is replaced by $f^s + f^l $ due to (\ref{eqn:PE_substit}) and 
(\ref{eqn:Dirac_splitting}) and $g$ is the Dirichlet boundary conditions defined in (\ref{eq:DH_solution}).

It was proved in  \cite{khor-DiracRS:2018}  and demonstrated in \cite{BeKKKS:18} that the function $f^s$ and the corresponding 
short-range potential $u^s$ are localized within the molecular region $\Omega_m$ and vanishes on the 
interface $\Gamma$. Moreover, from (\ref{eqn:diel_kappa_def}), the function $\bar{\kappa}$ is piecewise 
constant and $\bar{\kappa} = 0$ in $\Omega_m$. Therefore, we can rewrite the Boltzmann distribution term 
in (\ref{eqn:PBE_splitting}) as
\begin{equation}\label{eqn:Boltzmann_lr}
   \begin{aligned}
 & \bar{\kappa}^2\sinh(u^s + u^r) =   \bar{\kappa}^2\sinh(u^r), \quad  & \mbox{ because } \, u^s = 0 \,
  \mbox{ in } \Omega_s.
  \end{aligned}
  \end{equation}
Consequently, following the splitting of the Dirac delta distributions in (\ref{eqn:Dirac_splitting}), the 
short-range component of the potential satisfies the Poisson equation, i.e.,
\begin{equation}\label{eqn:PE_sr}
 -\nabla\cdot(\epsilon_m\nabla u^s) = f^s \quad \mbox{in} \,\, \mathbb{R}^3.
\end{equation}
It can be easily shown that \[u^s(\bar{x}) = {\bf P}_s \] is the cumulated canonical tensor in (\ref{eq:CCT_tensor}) 
which represents the precomputed short-range potential sum supported within the solute domain $\Omega_m$.

Subtracting (\ref{eqn:PE_sr}) from (\ref{eqn:PBE_splitting}) and using (\ref{eqn:Boltzmann_lr}), 
we obtain the NRPBE as follows
\begin{equation}\label{eqn:RPBE_decomp}
-\nabla \cdot(\epsilon\nabla u^r(\bar{x})) + \bar{\kappa}^2(\bar{x})\sinh(u^r(\bar{x})) 
   = f^l,   \quad \mbox{in} \,\,\, \Omega,
\end{equation}
subject to 
\begin{equation}\label{eqn:RPBE_bc}
u^r(\bar{x}) = \bar{\kappa}^2(\bar{x})\mathbf{P}_l   \quad \mbox{on} \,\, \partial{\Omega}.
\end{equation}

We recall that the regularization scheme for linear PBE introduced in \cite{BeKKKS:18}
reads as follows,
\begin{equation}\label{eqn:LRPBE_decomp}
-\nabla \cdot(\epsilon\nabla u^r(\bar{x})) + \bar{\kappa}^2(\bar{x})u^r(\bar{x}) 
   = f^l(\bar{x}),   \, \mbox{in} \,\,\,  \Omega,
\end{equation}
subject to the Dirichlet boundary conditions 
\begin{equation}\label{eqn:LRPBE_decomp_bc}
u^r(\bar{x}) = \bar{\kappa}^2(\bar{x})\mathbf{P}_l   \quad \mbox{on} \,\, \partial{\Omega}.
\end{equation}
In this way, (\ref{eqn:RPBE_decomp}) -- (\ref{eqn:RPBE_bc}) generalizes   the regularization scheme 
(\ref{eqn:LRPBE_decomp}) -- (\ref{eqn:LRPBE_decomp_bc}) to the nonlinear case.

Notice that by construction, the short-range potential vanishes on the interface $\Gamma$, hence it 
satisfies the discrete Poisson equation in (\ref{eqn:PE}) with the respective charge density $f^s$ and zero 
boundary conditions on $\Gamma$. Therefore, we recall (see \cite{BeKKKS:18} for the detailed discussion)  
that this equation can be subtracted from the full 
linear discrete PE system, such that the long-range component of the solution, $\mathbf{P}_l$, will satisfy 
the same linear system of equations (same interface conditions), but with a modified charge density 
corresponding to the weighted sum of the long-range tensors $f^l$ only.

\section{Numerical approach to solving the NRPBE} \label{sec:RPBE_numerical_approach}

Consider the uniform 3D $n^{\otimes 3}$ rectangular grid in $\Omega = [-b,b]^3$ with the mesh parameters 
$dx,dy,dz < 0.5$. One standard way of solving the NRPBE in (\ref{eqn:RPBE_decomp}) is that it is first 
discretized in space to obtain a nonlinear system in matrix-vector form 
\begin{equation}\label{eq:FOM}
 A(u_{\mathcal{N}}^r) = b^r,
 \qquad \mbox{ in }\, \mathbb{R}^{\mathcal{N}},
\end{equation}
where $A(u_{\mathcal{N}}^r) \in \mathbb{R}^{\mathcal{N}\times \mathcal{N}}$, $b^r 
\in \mathbb{R}^{\mathcal{N}}$, and $u_{\mathcal{N}}^r$ is the discretized solution vector. Here, 
$\mathcal{N}$ is usually in $\mathcal{O}(10^6)$.

Then system (\ref{eq:FOM}) can be solved using several existing techniques. For 
example, the nonlinear relaxation methods has been implemented in the Delphi software 
\cite{Rocchia_2001}, the nonlinear conjugate gradient (CG) method has been implemented in University 
of Houston Brownian Dynamics (UHBD) software \cite{Brock_1992}, the nonlinear multigrid (MG) method 
\cite{Oberoi_1993} and the inexact Newton method have been implemented in the adaptive 
Poisson-Boltzmann solver (APBS) software \cite{Holst:95}.  

In this study, we apply a different approach of solving (\ref{eqn:RPBE_decomp})
\cite{Mirzadeh:13, Shestakov:2002, Ji:2018}. In particular, an iterative approach is first applied to 
the continuous NRPBE in (\ref{eqn:RPBE_decomp}), where at the $(n+1)$st iteration step, the NRPBE is 
approximated by a linear equation via the Taylor series truncation. The expansion point of the Taylor 
series is the continuous solution $(u^r)^n$ at the $n$th iteration step.

Consider $(u^r)^n$ as the approximate solution at the $n$th iterative step, then the nonlinear 
term $\sinh((u^r)^{n+1})$ at the $(n+1)$st step is approximated by its truncated Taylor series 
expansion as follows
\begin{equation}\label{eq:Taylor_expansion_sinh}
\sinh((u^r)^{n+1}) \approx \sinh((u^r)^n) 
+ ((u^r)^{n+1} - (u^r)^n)\cosh((u^r)^n).
\end{equation}
Substituting the approximation (\ref{eq:Taylor_expansion_sinh}) into (\ref{eqn:RPBE_decomp}), we 
obtain
\begin{multline}\label{eq:PBE_approx_sinh}
-\nabla\cdot(\epsilon(\bar{x})\nabla (u^r)^{n+1}) 
+ \bar{\kappa}^2(\bar{x})\cosh((u^r)^n)(u^r)^{n+1} 
= -\bar{\kappa}^2(\bar{x})\sinh((u^r)^n) \\
+ \bar{\kappa}^2(\bar{x})\cosh((u^r)^n)(u^r)^n + b^r.
\end{multline}
The equation in (\ref{eq:PBE_approx_sinh}) is linear, and can then be numerically solved by first 
applying spatial discretization. In this regard, we first define
\begin{equation}\label{eqn:hyperbolic_cosine_vec}
\cosh\odot u_{\mathcal{N}}^r =: w = \begin{bmatrix}
		w_1 \\
		w_2 \\
		\vdots \\
		w_{\mathcal{N}}
		\end{bmatrix},
\end{equation}
where $\odot$ is the elementwise operation on a vector. 

Then, we construct the corresponding diagonal matrix from (\ref{eqn:hyperbolic_cosine_vec}) of the 
form \[B = \mbox{diag}(w_1, w_2, \ldots, w_{\mathcal{N}}).\] 
Finally, we obtain the following linear system 
\begin{equation}\label{eqn:Nonaffine_form_iterative_FOM}
 A_1(u_{\mathcal{N}}^r)^{n+1} + A_2B^{n}(u_{\mathcal{N}}^r)^{n+1} 
 = -A_2\sinh\odot(u_{\mathcal{N}}^r)^{n} + A_2B^{n}(u_{\mathcal{N}}^r)^{n} + b_1^r + b_2,
\end{equation}
where $A_1$ is the Laplacian matrix and $A_2$ is a diagonal matrix containing the $\bar{\kappa}^2$ 
function. Note that the diagonal matrix $B^{n}$ changes at each iteration step, therefore, it cannot 
be precomputed. The vectors $b_1^r$ and $b_2$ are the regularized approximation of the Dirac delta 
distributions and the Dirichlet boundary conditions, respectively.

Let 
\begin{equation}\label{eqn:affine_A_iter}
 A(\cdot) = A_1 + A_2B^n 
\end{equation} 
and
\begin{equation}\label{eqn:affine_F_iter}
 F : \mbox{right-hand side of} \: (\ref{eqn:Nonaffine_form_iterative_FOM}).
\end{equation} 
Then we obtain
\begin{equation}\label{eq:PBE_system_iterative}
A((u_{\mathcal{N}}^r)^n)(u_{\mathcal{N}}^r)^{n+1} = F((u_{\mathcal{N}}^r)^n), 
\quad n = 0,1, \ldots.
\end{equation} 

Then, at each iteration, system (\ref{eq:PBE_system_iterative}) is a linear system w.r.t. 
$(u^r_{\mathcal{N}})^{n+1}$, which can be solved by any linear system solver of choice. In this study, 
we employ the aggregation-based algebraic multigrid method (AGMG) \footnote{\textbf{AGMG} implements an 
aggregation-based algebraic multigrid method, which solves algebraic systems of linear equations, and is 
expected to be efficient for large systems arising from the discretization of scalar second order 
elliptic PDEs \cite{Notay:2010}.} \cite{Notay:2010}. \Cref{alg:Iterative_NRPBE} summarizes the detailed 
iterative approach of solving (\ref{eq:PBE_system_iterative}). This approach of first linearization, 
then discretization is shown to be more efficient than the standard way of first discretization and 
then linearization, via, for example, the Newton iteration. The advantage of the proposed approach is 
that it avoids computing the Jacobian of a huge matrix. It is observed that it converges faster than 
the standard Newton approach. 

\begin{algorithm}[t]
  \caption{Iterative solver for the NRPBE}\label{alg:Iterative_NRPBE}
    \begin{algorithmic}[1]
    \Require Initialize the potential $(u_{\mathcal{N}}^r)^0$, e.g., $(u_{\mathcal{N}}^r)^0 = 0$ 
	      and the tolerance $\delta^0 = 1$. 
    \Ensure The converged NRPBE solution $(u_{\mathcal{N}}^r)^n$ at $\delta^n \leq \tau$.
    \While{$\delta^n \geq \tau$}
    \State Solve the linear system (\ref{eq:PBE_system_iterative}) for 
      $(u_{\mathcal{N}}^r)^{n+1}$ using AGMG.
    \State $\delta^{n+1} \gets \|(u_{\mathcal{N}}^r)^{n+1} - (u_{\mathcal{N}}^r)^n\|_2$.
    \State $(u_{\mathcal{N}}^r)^n \gets (u_{\mathcal{N}}^r)^{n+1}$.
    \EndWhile
    \State \textbf{end while}
  \end{algorithmic}
\end{algorithm}

The benefits of the RS tensor format as a solution decomposition technique over the existing techniques 
in the literature are highlighted as follows. First, the efficient splitting of the short- and long-range 
parts in the target tensor circumvents the need to modify jump conditions at the interface and the 
use of $\epsilon$ and $\bar{\kappa}$ as cut-off functions, e.g., in (\ref{eqn:PBE_solndecomp}). Secondly, 
the long-range part in the RS tensor decomposition of the Dirac delta distributions 
\cite{khor-DiracRS:2018} vanishes at the interface and, therefore, the modified charge density in 
(\ref{eqn:Dirac_splitting}) generated by this long-range component remains localized in the solute 
region. Thirdly, the boundary conditions are obtained from $\mathbf{P}_l$, the long-range part of the 
free space potential sum, thereby avoiding the computational costs involved in solving some 
external analytical function at the boundary. Lastly, only a single system of algebraic equations is 
solved for the regularized component of the collective potential which is then added to the directly 
precomputed short-range contribution, $u^s(\bar{x})$. This is more efficient than, for instance, in 
\cite{Xie:14}, where the regularized PBE model is subdivided into the linear interface and the 
nonlinear interface problems which are solved independently, with respective boundary and interface 
conditions.

\section{Numerical results}
\label{sec:Numer_Tests}

In this section, we consider $n^{\otimes 3}$ 3D uniform Cartesian grids  in a box $[-b,b]^3$ with equal step size 
$h =2b/(n-1)$ for computing the electrostatic potentials of the PBE on a modest PC with the following 
specifications: Intel (R) Core (TM) $i7-4790$ CPU @ 3.60GHz with 8GB RAM. The FDM is used to discretize the 
PBE in this work and the numerical computations are implemented in the MATLAB software, version R2017b.

\subsection{Numerical results for LPBE}
\label{ssec:Numer_Tests}

\begin{figure}[t]
\captionsetup{width=\linewidth}
    \centering
    \begin{subfigure}[b]{0.42\textwidth}
        \centering
        \includegraphics[width=\textwidth]{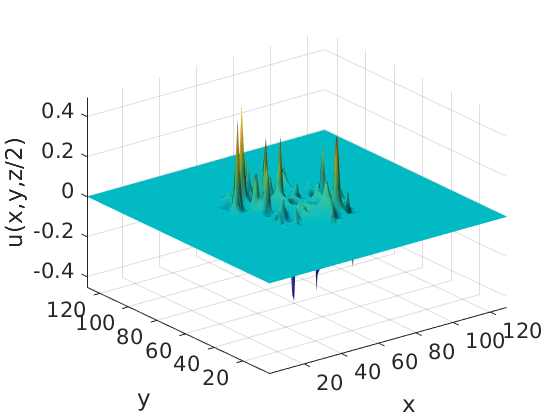}
        \caption{LPBE solution by the FDM solver.}
        \label{fig:LPBE_FDM}
    \end{subfigure}
    \hfill
    \begin{subfigure}[b]{0.42\textwidth}
        \centering
        \includegraphics[width=\textwidth]{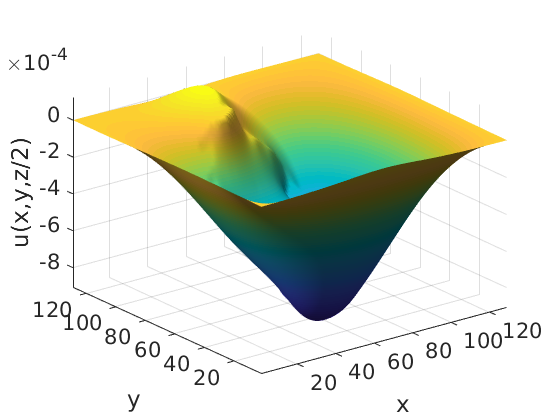}
        \caption{APBS vs FDM solution error.}
        \label{fig:Error_APBS_FDM}
    \end{subfigure}
    \caption{\label{fig:APBS_FDM_LPBE}
    The electrostatic potential for the protein Fasciculin 1 computed by the FDM solver (left) 
and the error between the APBS and FDM solutions (right) on $129^{\otimes 3}$ grid, at $0.15M$ ionic 
strength.}   
\end{figure}

First, we validate our FDM solver for the classical LPBE by comparing its solution with that of 
the APBS software package (version 1.5-linux64), which uses the multigrid (PMG) accelerated FDM 
\cite{Bakersept2001}. Here, we consider the protein Fasciculin 1, with 1228 atoms.
\Cref{fig:APBS_FDM_LPBE} shows the electrostatic potential of the PBE on a $n\times n$ grid surface 
with $n=129$ at the cross-section of the volume box ($60\,\mbox{\AA}$) in the middle of the $z$-axis computed 
by the FDM solver and the corresponding error between the two solutions. Here, we use the ionic strength 
of $0.15M$ and the dielectric coefficients $\epsilon_m =2$ and $\epsilon_s = 78.54$, respectively. The 
numerical results show that the FDM solver provides as accurate results as those of the APBS with a 
discrete $L_2$ error of $\mathcal{O}(10^{-4})$ in the full solution. 

The corresponding electrostatic potential energy for the aforementioned LPBE solvers on a sequence of 
fine grids is given in the \Cref{table:Energy_compare}. The results for solvation free energy of protein 
varieties are presented in \cite{KwBFM2017}. To validate the claim in \Cref{remk:LPBE_NPBE_diff}, we 
provide in the \Cref{table:Energy_LPBE_NPBE}, the comparison between the total electrostatic 
potential energies $\Delta G_{elec}$ in kJ/mol, between the LPBE and the nonlinear PBE (NPBE) 
computations on a sequence of fine grids using the APBS software package.
\begin{table}[t]
\centering
\captionsetup{width=\linewidth}
\begin{tabular}{|c|c|c|c|c|}
  \hline
  $h$& $\mathcal{N}$ & $\Delta G_{elec}$, FDM & $\Delta G_{elec}$, APBS 
  & Relative error\\ \hline
  {0.465} & $129^3$ & 91,232.9217 & 91,228.0388 & 5.3524e-5\\ \hline
  {0.375} & $161^3$ & 130,611.0021 & 130,606.0444 & 3.7962e-5\\ \hline
  {0.320} & $193^3$ & 170,159.4204 & 170,154.3821 & 2.9610e-5\\ \hline
\end{tabular}
\caption{Comparison of the total electrostatic potential energies $\Delta G_{elec}$ in kJ/mol, between FDM 
and APBS on a sequence of fine grids.}
\label{table:Energy_compare}
\end{table}

\begin{remark}\label{remk:LPBE_NPBE_diff}
 We reiterate that the solutions obtained from the LPBE and the nonlinear PBE are very close to each 
 other, even when the linearization condition does not hold \cite{Fogolari2002}. This is especially 
 manifested in protein molecules whose charge densities are small. However, in biomolecules with large 
 charge densities, for example, the DNA, significant differences might be observed at the solute-solvent 
 interface \cite{Fogolari99, Fogolari2002}. Moreover, the solution of the LPBE is usually used as the 
 initial guess for the nonlinear PBE.
\end{remark}
\begin{table}[t]
\centering
\captionsetup{width=\linewidth}
\begin{tabular}{|c|c|c|c|c|}
  \hline
  $h$& $\mathcal{N}$ & $\Delta G_{elec}$, LPBE & $\Delta G_{elec}$, NPBE 
  & Relative error\\ \hline
  {0.465} & $129^3$ & 91,228.0575 & 91,227.8354 & 2.4345e-6\\ \hline
  {0.375} & $161^3$ & 130,606.0630 & 130,605.8448 & 1.6707e-6\\ \hline
  {0.320} & $193^3$ & 170,154.4401 & 170,154.1862 & 1.4922e-6\\ \hline
\end{tabular}
\caption{Comparison of the total electrostatic potential energies $\Delta G_{elec}$ in kJ/mol, between the 
LPBE and the NPBE computations on a sequence of fine grids.}
\label{table:Energy_LPBE_NPBE}
\end{table}

\begin{remark}\label{remk:Energy_gridsize}
 Notice from \Cref{table:Energy_compare} and \Cref{table:Energy_LPBE_NPBE} that the electrostatic 
 potential energies $\Delta G_{elec}$ increase with decreasing grid/mesh size, $h$. This is caused by the 
 short-range electrostatic potential behaviour in $1/\|\bar{x}\|$ as $\|\bar{x}\|$ $ \to 0$. 
\end{remark}

\subsection{Accuracy of the nonlinear RPBE based on the RS tensor format}
\label{ssec:RPBE_RS_accuracy} 

Here, we provide the results for the calculation of electrostatic potential for the nonlinear RPBE (NRPBE) 
based on the RS tensor format and compare the results with those of the traditional NPBE for various 
proteins. First, we consider the protein Fasciculin 1 consisting of 1228 atoms of varying atomic radii as shown in 
\Cref{table:Fasciculin_atoms}. Notice that 322 of the total atoms have zero radius, which implies that 
we must annihilate them from the RS tensor format calculations so that they are not assigned Newton kernels. 
Therefore, we consider the smallest atom in the protein as that with 1 $\mbox{\AA}$ radius, (i.e., the 
Hydrogen atom).
\begin{table}[t]
\centering
\captionsetup{width=\linewidth}
\begin{tabular}{|c|c|c|c|c|c|c|c|}
  \hline
  \multicolumn{8}{|c|}{Atomic radii in $\mbox{\AA}$} \\ \hline
   Atomic radius & 0.00 & 1.00 & 1.40 & 1.50 & 1.70 & 1.85 & 2.00\\ \hline
  Number of atoms & 322 & 333 & 195 & 82 & 104 & 10 & 182\\ \hline
\end{tabular}
\caption{Atomic radii and the corresponding number of atoms for the constituent atoms of the protein 
Fasciculin 1.}
\label{table:Fasciculin_atoms}
\end{table}

\begin{figure}[t]
\centering
\captionsetup{width=\linewidth}
\resizebox{16cm}{10.0cm}{%
\begin{tikzpicture}

\node (APBS_protein_pot_129) at (0,3.0) {\begin{tabular}{l}
  {\includegraphics[width=0.48\textwidth]{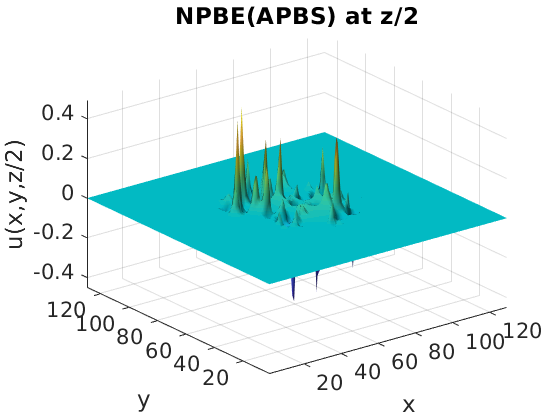}}
 \end{tabular}
 };
 
 \node (NRPBE_protein_pot_129) at (0,-3.0) {\begin{tabular}{l}
  {\includegraphics[width=0.48\textwidth]{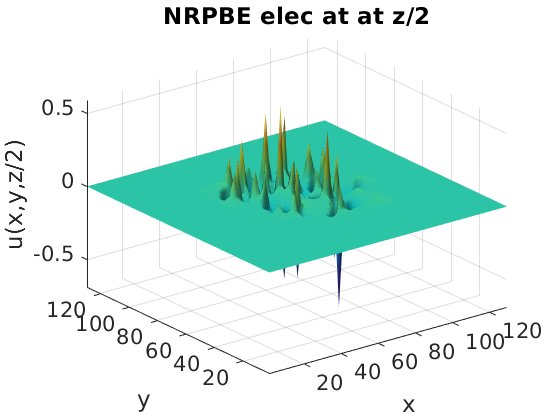}}
 \end{tabular}
 };
 
 \node (Error_NPBE_NRPBE_129) at (10,0) {\begin{tabular}{l}
  {\includegraphics[width=0.48\textwidth]{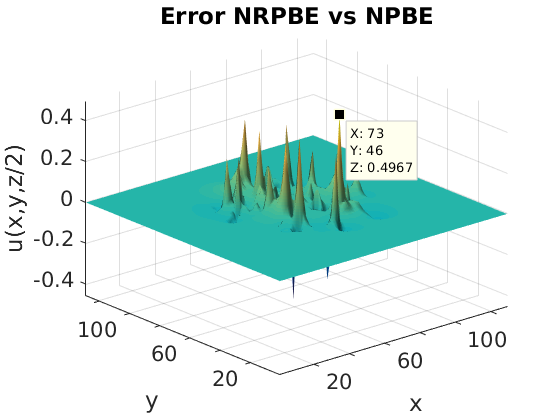}}
 \end{tabular}
 };
 \draw[->,draw=blue,thick] (3.5,0) -- (Error_NPBE_NRPBE_129.west);
 
 \end{tikzpicture}}
 \caption{\label{fig:APBS_NPBE_NRPBE_solns}
 Absolute error between the solutions of the traditional NPBE and the NRPBE for the protein 
 Fasciculin 1.}
\end{figure}

We provide the comparisons between the electrostatic potential computed by the NRPBE, based on the RS tensor format, 
with that of the traditional NPBE. \Cref{fig:APBS_NPBE_NRPBE_solns} shows the solutions from the 
two models and the corresponding error on $129^{\otimes 3}$ uniform Cartesian grid and a 
$60\mbox{\AA}$ domain length, at $0.15M$ ionic strength.
\begin{remark}\label{rmk:molecular_solute_solns}
Notice that the error is predominant within the molecular region, where the solution is singular. 
However, in the solute region, which is dominated by the long-range regime, the error is small, of 
order $\mathcal{O}(10^{-5})$. 
\end{remark}

\begin{figure}[t]
\captionsetup{width=\linewidth}
    \centering
    \begin{subfigure}[b]{0.42\textwidth}
        \centering
        \includegraphics[width=\textwidth]{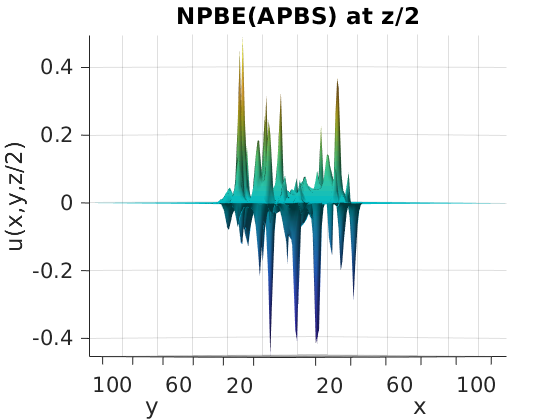}
        \caption{Cross-sectional view of NPBE solution in \Cref{fig:APBS_NPBE_NRPBE_solns}.}
        \label{fig:NPBE_APBS}
    \end{subfigure}
    \hfill
    \begin{subfigure}[b]{0.42\textwidth}
        \centering
        \includegraphics[width=\textwidth]{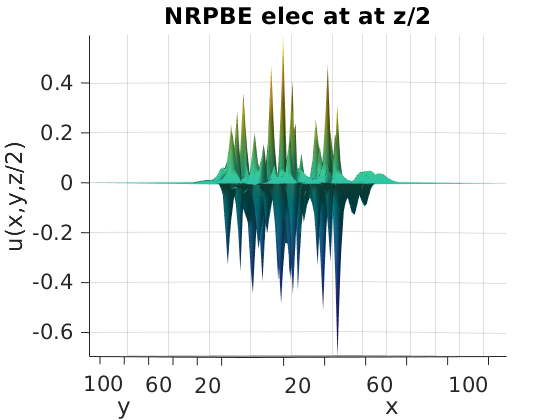}
        \caption{Cross-sectional view of NRPBE solution in \Cref{fig:APBS_NPBE_NRPBE_solns}.}
        \label{fig:NRPBE_FDM}
    \end{subfigure}
    \caption{\label{fig:APBS_FDM_NPBE_view}
    The cross-sectional view of the electrostatic potentials in \Cref{fig:APBS_NPBE_NRPBE_solns}.}   
\end{figure}

\Cref{fig:APBS_FDM_NPBE_view} provides the cross-sectional view of the electrostatic potential shown 
in \Cref{fig:APBS_NPBE_NRPBE_solns}, for demonstrating the accuracy of the numerical treatment of the 
solution singularities inherent in the NRPBE model as compared with 
the traditional NPBE model. Notice that the NRPBE is capable of capturing 
exactly, the short-range component of the total potential sum because this part is precomputed analytically 
thereby avoiding the numerical errors generated by the traditional NPBE solver.

\begin{remark}\label{rmk:singularities_accuracy}
\Cref{fig:NRPBE_FDM} contains densely populated singularities/cusps as a result of explicit 
treatment of each atomic charge by the short-range part of the RS tensor 
whereas \Cref{fig:NPBE_APBS}, displays sparsely populated singularities, most of which are not sharp due 
to the redundant smoothing/smearing effect of the atomic charges by the cubic spline interpolation.
\end{remark}

\begin{figure}[H]
\centering
\captionsetup{width=\linewidth}
\resizebox{16cm}{10.0cm}{%
\begin{tikzpicture}

\node (APBS_protein_pot_129) at (0,3.0) {\begin{tabular}{l}
  {\includegraphics[width=0.48\textwidth]{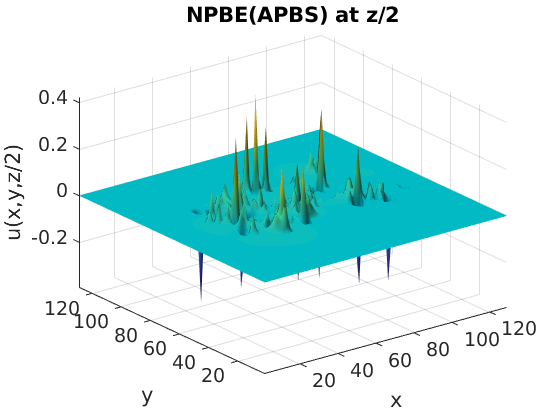}}
 \end{tabular}
 };
 
 \node (NRPBE_protein_pot_129) at (0,-3.0) {\begin{tabular}{l}
  {\includegraphics[width=0.48\textwidth]{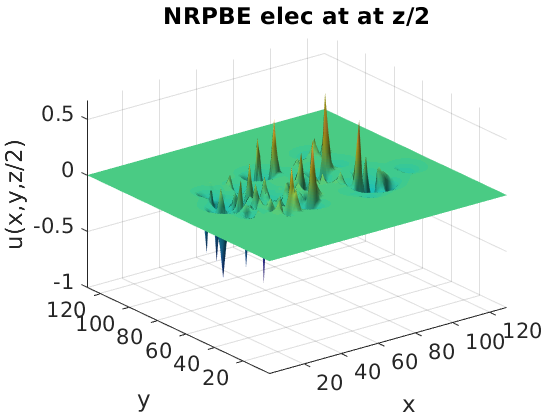}}
 \end{tabular}
 };
 
 \node (Error_NPBE_NRPBE_129) at (10,0) {\begin{tabular}{l}
  {\includegraphics[width=0.48\textwidth]{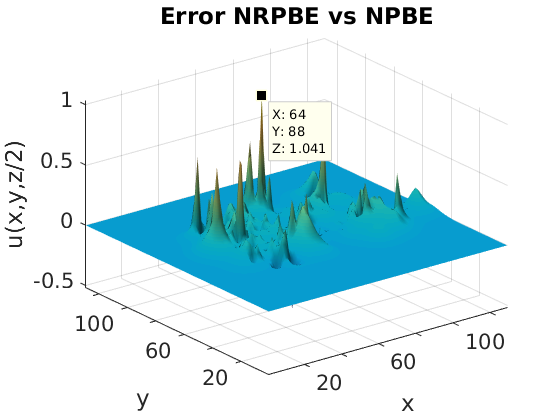}}
 \end{tabular}
 };
 \draw[->,draw=blue,thick] (3.5,0) -- (Error_NPBE_NRPBE_129.west);
 
 \end{tikzpicture}}
 \caption{\label{fig:APBS_NPBE_NRPBE_solns_2809}
 Absolute error between the solutions of the traditional NPBE and the NRPBE for the 
 murine-human chimera of leukemia inhibitory factor (LIF).}
\end{figure}

Secondly, we provide results for a 180-residue cytokine solution NMR structure of a murine-human chimera 
of leukemia inhibitory factor (LIF) \cite{HinMauZhaNic:98} consisting of 2809 atoms. The corresponding 
variation in atomic radii and the corresponding atomic occurrences are shown in \Cref{table:2809_atoms}.  
\Cref{fig:APBS_NPBE_NRPBE_solns_2809} shows the comparison between the electrostatic potential of NRPBE, 
with that of the classical NPBE and the corresponding error on a $129^{\otimes 3}$ grid and a 
$65\mbox{\AA}$ domain length, at $0.15M$ ionic strength. 
\begin{table}[t]
\centering
\captionsetup{width=\linewidth}
\begin{tabular}{|c|c|c|c|c|c|c|c|}
  \hline
  \multicolumn{8}{|c|}{Atomic radii in $\mbox{\AA}$} \\ \hline
   Atomic radius & 0.2245 & 0.4500 & 0.9000 & 1.3200 & 1.3582 & 1.4680 &$\geq$ 1.7000\\ \hline
  Number of atoms & 315 & 6 & 6 & 1032 & 54 & 6 & 1390\\ \hline
\end{tabular}
\caption{Atomic radii and the corresponding number of atoms for the constituent atoms of a 180-residue 
cytokine solution NMR structure of a murine-human chimera of leukemia inhibitory factor (LIF).}
\label{table:2809_atoms}
\end{table}

\begin{figure}[H]
\captionsetup{width=\linewidth}
    \centering
    \begin{subfigure}[b]{0.42\textwidth}
        \centering
        \includegraphics[width=\textwidth]{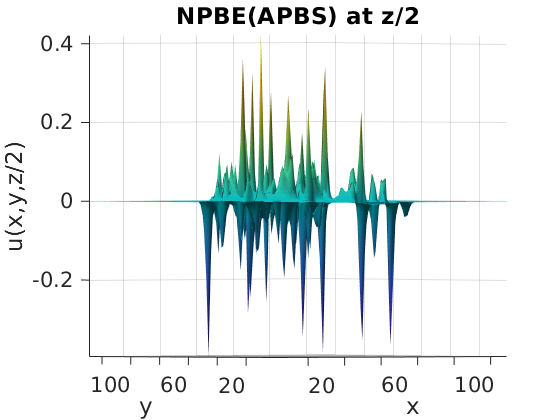}
        \caption{Cross-sectional view of NPBE solution in the \Cref{fig:APBS_NPBE_NRPBE_solns}.}
        \label{fig:NPBE_APBS_2809}
    \end{subfigure}
    \hfill
    \begin{subfigure}[b]{0.42\textwidth}
        \centering
       \includegraphics[width=\textwidth]{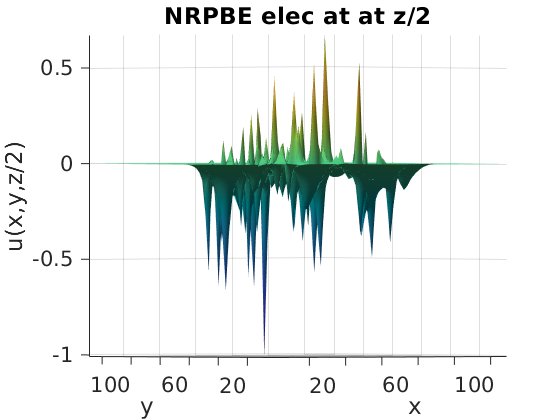}
        \caption{Cross-sectional view of NRPBE solution in the \Cref{fig:APBS_NPBE_NRPBE_solns}.}
        \label{fig:NRPBE_FDM_2809}
    \end{subfigure}
    \caption{\label{fig:APBS_FDM_NPBE_view_2809}
    The cross-sectional view of the electrostatic potentials in the 
    \Cref{fig:APBS_NPBE_NRPBE_solns_2809}.}   
\end{figure}

\begin{remark}\label{rmk:accuracy_2809}
In a similar vein, we notice in \Cref{fig:APBS_FDM_NPBE_view_2809} that the error is predominant 
within the molecular region, where the solution is singular. It is also worth mentioning that the small atomic radii, ($< 0.9 \mbox{\AA}$), 
in \Cref{table:2809_atoms} are treated independently in terms of the RS tensor splitting of the 
short- and long-range potentials.
\end{remark}

\subsection{Runtimes and Computational Speed-ups}
\label{ssec:RPBE_RS_runtimes}

We compare the runtimes of computing both the classical and regularized PBE models in 
\Cref{table:Runtimes_speedups} for the protein Fasciculin 1 in an $n^3=129^3$ domain of $60\,\mbox{\AA}$ 
length at an ionic strength of $0.15M$. Notice that the runtimes for the LPBE and the LRPBE are almost 
equal because the linear systems are solved by the same solver (i.e., AGMG). On the other hand, the 
runtime for solving the nonlinear system for the NRPBE is half that of the NPBE due to the absence of 
the Dirac delta distributions and their corresponding solution singularities in 
our scheme, and  which increase  the computational costs in NPBE.  

\begin{table}[t]
\centering
\captionsetup{width=\linewidth}
\begin{tabular}{|c|c|c|c|}
  \hline
  \multicolumn{4}{|c|}{Runtime (seconds) and speed-up} \\
  \hline
   & LPBE & LRPBE & Speed-up\\
  \hline
  Solve linear system & 5.26 & 6.34 & $\approx$ 1\\ \hline
  Total runtime & 15.25 & 16.47 & $\approx$ 1\\ \hline
   & NPBE & NRPBE & Speed-up\\
  \hline
   Solve nonlinear system & 24.23 & 12.30 & 1.97\\ \hline
  Total runtime & 34.40 & 28.30 & 1.21\\ \hline
\end{tabular}
\caption{Runtimes and speed-ups for LPBE, LRPBE, NPBE, and NRPBE.}
\label{table:Runtimes_speedups}
\end{table}

\section{Conclusions}\label{sec:Conclusions}

In this paper, we apply the RS tensor format for a solution decomposition 
of the nonlinear PBE for computation of electrostatic potential of large solvated biomolecules. 
The efficacy of the tensor-based regularization scheme established in \cite{BeKKKS:18} for the linear PBE,  
is based on the unprecedented properties of the grid-based RS tensor splitting 
of the Dirac delta distribution \cite{khor-DiracRS:2018}. Similar to the linear case, the key computational 
benefits are attributed to the localization of the modified Dirac delta distributions within the molecular 
region and the automatic maintaining of the continuity of the Cauchy data on the solute-solvent interface. 
Moreover, our computational scheme entails solving only a single system of algebraic equations for the 
regularized component of the collective electrostatic potential discretized by the FDM.
The total potential is obtained by adding this solution to the directly 
precomputed low-rank tensor representation of the short-range contribution. 
 
The main properties of the presented scheme are demonstrated by various numerical tests. 
For instance, \Cref{fig:APBS_FDM_NPBE_view} and \Cref{fig:APBS_FDM_NPBE_view_2809} 
vividly demonstrate that the traditional PBE model does not accurately capture the solution 
singularities which originate from the short-range component of the total target electrostatic 
potential in the numerical approximation. In the RPBE, the Dirac delta distribution is replaced by a 
smooth long-range function from (\ref{eqn:Dirac_splitting}). It only requires one to solve for the 
long-range electrostatic potential numerically and add this solution to the short-range component 
which is computed a priori using the canonical tensor approximation to the Newton kernel. The 
resultant total potential sum is of high accuracy as demonstrated by \Cref{fig:NRPBE_FDM} and 
\Cref{fig:NRPBE_FDM_2809}.

\section*{Acknowledgement}
The authors thank the following organizations for financial and material support on 
this project: International Max Planck Research School (IMPRS) for Advanced Methods in Process 
and Systems Engineering and Max Planck Society for the Advancement of Science (MPG).

  \begin{footnotesize}

    \bibliographystyle{unsrtnat}
  \bibliography{BSE_Fock_Sums1.bib,Kweyu_refer.bib}

  \end{footnotesize}

\end{document}